\documentclass[11pt]{article}
\usepackage[centertags]{amsmath}
\usepackage{amsfonts}
\usepackage{amssymb}
\usepackage{amsthm}
\usepackage{newlfont}
\usepackage{graphicx}

\newfont{\bb}{msbm10 at 12pt}
\def\r{\hbox{\bb R}}

\def\e{\hbox{\bf E}}
\def\t{\hbox{\bf T}}
\def\n{\hbox{\bf N}}
\def\b{\hbox{\bf B}}

\newtheorem{theorem}{Theorem}[section]

\newtheorem{remark}[theorem]{Remark}

\newtheorem{lemma}[theorem]{Lemma}
\newtheorem{example}[theorem]{Example}

\begin{document}

\title{Position vectors of slant helices in Euclidean space E$^3$}
\author{ Ahmad T. Ali\\Mathematics Department\\
 Faculty of Science, Al-Azhar University\\
 Nasr City, 11448, Cairo, Egypt\\
email: atali71@yahoo.com}

\maketitle
\begin{abstract}
In classical differential geometry, the problem of the determination of the position vector of an arbitrary space curve according to the intrinsic equations $\kappa=\kappa(s)$ and $\tau=\tau(s)$ (where $\kappa$ and $\tau$ are the curvature and torsion of the space curve $\psi$, respectively) is still open \cite{eisenh, lips}. However, in the case of a plane curve, helix and general helix, this problem is solved. In this paper, we solved this problem in the case of a slant helix. Also, we applied this method to find the representation of a Salkowski, anti-Salkowski curves and a curve of constant precession, as examples of a slant helices, by means of intrinsic equations.
\end{abstract}

\emph{MSC:} 53C40, 53C50

\emph{Keywords}: Classical differential geometry; Frenet equations;  Slant helices; Intrinsic equations.

\section{Introduction }
In the local differential geometry, we think of curves as a geometric set of points, or locus. Intuitively, we are thinking of a curve as the path traced out by a particle moving in $\e^3$. So, the investigating position vectors of the curves in a classical aim to determine behavior of the particle (curve).

Helix is one of the most fascinating curves in science and nature. Scientist have long held a fascinating, sometimes bordering on mystical obsession, for helical structures in nature. Helices arise in nano-springs, carbon nano-tubes, $\alpha$-helices, DNA double and collagen triple helix, lipid bilayers, bacterial flagella in salmonella and escherichia coli, aerial hyphae in actinomycetes, bacterial shape in spirochetes, horns, tendrils, vines, screws, springs, helical staircases and sea shells \cite{choua, lucas, watson}. Also we can see the helix curve or helical structures in fractal geometry, for instance hyperhelices \cite{toledo}. In the field of computer aided design and computer graphics, helices can be used for the tool path description, the simulation of kinematic motion or the design of highways, etc. \cite{yang}. From the view of differential geometry, a helix is a geometric curve with non-vanishing constant curvature $\kappa$ and non-vanishing constant torsion $\tau$ \cite{barros}. The helix may be called a {\it circular helix} or {\it W-curve} \cite{ilarslan, mont1}.

Indeed a helix is a special case of the {\it general helix}. A curve of constant slope or general helix in Euclidean 3-space $\e^3$ is defined by the property that the tangent makes a constant angle with a fixed straight line called the axis of the general helix. A classical result stated by Lancret in 1802 and first proved by de Saint Venant in 1845 (see \cite{struik} for details) says that: {\it A necessary and sufficient condition that a curve be a general helix is that the ratio $$\dfrac{\kappa}{\tau}$$ is constant along the curve, where $\kappa$ and $\tau$ denote the curvature and the torsion, respectively}. A general helices or {\it inclined curves} are well known curves in classical differential geometry of space curves \cite{milm} and we refer to the reader for recent
works on this type of curves \cite{ali2, barros, gl1, mont2, sc, tur2}.

Izumiya and Takeuchi \cite{izumi} have introduced the concept of slant helix by saying that the normal lines make a constant angle with a fixed straight line. They characterize a slant helix if and only if the geodesic curvature of the principal image of the principal normal indicatrix
$$
\frac{\kappa^2}{(\kappa^2+\tau^2)^{3/2}}\Big(\frac{\tau}{\kappa}\Big)'
$$
is a constant function. They also called a curve is canonical geodesic if
$$
\Big(\frac{\tau}{\kappa}\Big)'
$$
is a constant function. Kula and Yayli \cite{kula} have studied spherical images of tangent indicatrix and binormal indicatrix of a slant helix and they showed that the spherical images are spherical helices. Recently, Kula et al. \cite{kula1} investigated the relations between a general helix and a slant helix. Moreover, they obtained some differential equations which they are characterizations for a space curve to be a slant helix.

Many important results in the theory of the curves in $\e^3$ were initiated by G. Monge and G. Darboux pioneered the moving frame idea. Thereafter, F. Frenet defined his moving frame and his special equations which play important role in mechanics and kinematics as well as in differential geometry \cite{boyer}.

For a unit speed curve with non-vanishing curvature $\kappa \neq 0$, it is well-known the following result \cite{hacis}:

\begin{theorem}\label{th-main} A curve is defined uniquely by its curvature and torsion as function of a natural parameters.
\end{theorem}

The equations
$$
\kappa=\kappa(s),\,\,\,\,\,\tau=\tau(s)
$$
which give the curvature and torsion of a curve as functions of $s$ are called the {\it natural} or {\it intrinsic equations} of a curve, for they completely define the curve.

Given two functions of one parameter (potentially curvature and torsion parameterized by arc-length) one might like to find an arc-length parameterized curve for which the two functions work as the curvature and the torsion. This problem, known as {\it solving natural equations}, is generally achieved by solving a {\it Riccati equation} \cite{struik}. Barros et. al. \cite{ba2} showed that the general helices in Euclidean 3-space $\e^3$ and in the three-sphere $\bold{S}^3$ are geodesic either of right cylinders or of Hopf cylinders according to whether the curve lies in $\e^3$ or $\bold{S}^3$, respectively.

In classical differential geometry, the problem of the determination of parametric representation of the position vector of an arbitrary space curve according to the intrinsic equations is still open \cite{eisenh, lips}. This problem is solved in three special cases only. Firstly, in the case of a plane curve $(\tau=0)$ \cite{struik}. Secondly, in the case of a circular helix ($\kappa$ and $\tau$ are both non-vanishing constants) \cite{struik}. Recently, Ali \cite{ali1} solved this problem in the case of a general helix . However, this problem is not solved in other cases of the space curve or in the general case.

Our main result in this work is to determine the parametric representation of the  position vector $\psi$ from intrinsic equations in $\e^3$ for the case of a slant helix.

\section{Preliminaries }
In Euclidean space $\e^3$, it is well known that to each unit speed curve with at least four continuous derivatives, one can associate three mutually orthogonal unit vector fields $\t$, $\n$ and $\b$ are respectively, the tangent, the principal normal and the binormal vector fields \cite{hacis}.

We consider the usual metric in Euclidean 3-space $\e^3$, that is,
$$
\langle,\rangle=dx_1^2+dx_2^2+dx_3^2,
$$
where $(x_1,x_2,x_3)$ is a rectangular coordinate system of $\e^3$.  Let $\psi:I\subset\r\rightarrow\e^3$, $\psi=\psi(s)$, be an arbitrary curve in $\e^3$. The curve $\psi$ is said to be of unit speed (or parameterized by the  arc-length) if $\langle\psi'(s),\psi'(s)\rangle=1$ for any $s\in I$. In particular, if $\psi(s)\not=0$ for any $s$, then it is possible to re-parameterize $\psi$, that is, $\alpha=\psi(\phi(s))$ so that $\alpha$ is parameterized by the arc-length. Thus, we will assume throughout this work that $\psi$ is a unit speed curve.

Let $\{\t(s),\n(s),\b(s)\}$ be the moving frame along $\psi$, where the vectors $\t, \n$ and $\b$ are mutually orthogonal vectors satisfying $\langle\t,\t\rangle=\langle\n,\n\rangle=\langle\b,\b\rangle=1$.
The Frenet equations for $\psi$ are given by (\cite{struik})
\begin{equation}\label{u1}
 \left[
   \begin{array}{c}
     \t'(s) \\
     \n'(s) \\
     \b'(s) \\
   \end{array}
 \right]=\left[
           \begin{array}{ccc}
             0 & \kappa(s) & 0 \\
             -\kappa(s) & 0 & \tau(s) \\
             0 & -\tau(s) & 0 \\
           \end{array}
         \right]\left[
   \begin{array}{c}
     \t(s) \\
     \n(s) \\
     \b(s) \\
   \end{array}
 \right].
 \end{equation}

If $\tau(s)=0$ for any $s\in I$, then $\b(s)$ is a constant vector $V$ and the curve $\psi$ lies in a $2$-dimensional affine subspace orthogonal to $V$, which is isometric to the Euclidean $2$-space $\e^{2}$.

We observe that the Frenet equations form a system of three vector differential equations of the first order in $\t, \n$ and $\b$. It is reasonable to ask, therefore, given arbitrary continuous functions $\kappa$ and $\tau$, whether or not there exist solutions $\t, \n, \b$ of the Frenet equations, and hence, since $\psi'=\t$, a curve
\begin{equation}\label{u111}
\begin{array}{ll}
\psi(s)=\int\,\t(s)\,ds+\bold{C},
  \end{array}
\end{equation}
which the prescribed curvature and torsion. The answer is in the affirmative and is given by:

\begin{theorem}{\bf (Fundamental existence and uniqueness theorem for space curve).}\label{th-main} Let $\kappa(s)$ and  $\tau(s)$ be arbitrary continuous function on $a\leq s \leq b$. Then there exists, except for position in space, one and only one curve $C$ for which $\kappa(s)$ is the curvature, $\tau(s)$ is the torsion and $s$ is a natural parameter along $C$.
\end{theorem}

\section{Position vectors of a space curves }

\begin{theorem}\label{th-main} Let $\psi=\psi(s)$ be an unit speed curve parameterized by the arclength $s$. Suppose $\psi=\psi(\theta)$ is another parametric representation of this curve by the parameter $\theta=\int\kappa(s)ds$.  Then, the principal normal vector $\n$ satisfies a vector differential equation of third order as follows:
\begin{equation}\label{u2}
\frac{1}{f(\theta)}\Bigg[
\frac{1}{f'(\theta)}\Bigg(\n''(\theta)+\Big(1+f^2(\theta)\Big)\n(\theta)\Bigg)\Bigg]'+\n(\theta)=0,
\end{equation}
where $f(\theta)=\frac{\tau(\theta)}{\kappa(\theta)}$.
\end{theorem}

{\bf Proof.} Let $\psi=\psi(s)$ be an unit speed curve. If we write this curve in the another parametric representation $\psi=\psi(\theta)$, where $\theta=\int\kappa(s)ds$, we have new Frenet equations as follows:
\begin{equation}\label{u3}
 \left[
   \begin{array}{c}
     \t'(\theta) \\
     \n'(\theta) \\
     \b'(\theta) \\
   \end{array}
 \right]=\left[
           \begin{array}{ccc}
             0 & 1 & 0 \\
             -1 & 0 & f(\theta) \\
             0 & -f(\theta) & 0 \\
           \end{array}
         \right]\left[
   \begin{array}{c}
     \t(\theta) \\
     \n(\theta) \\
     \b(\theta) \\
   \end{array}
 \right],
 \end{equation}
where $f(\theta)=\frac{\tau(\theta)}{\kappa(\theta)}$. If we differentiate the second equation of the new Frenet equations (\ref{u3}) and using the first and the third equations, we have
\begin{equation}\label{u4}
\b(\theta)=\dfrac{1}{f'(\theta)}\Big[\n''(\theta)+\Big(1+f^2(\theta)\Big)\n(\theta)\Big].
\end{equation}
Differentiating the above equation and using the last equation from (\ref{u3}), we obtain a vector differential equation of third order (\ref{u2}) as desired.

The equation (\ref{u2}) is not easy to solve in general case. If one solves this equation, the natural representation of the position vector of an arbitrary space curve can be determined as follows:
\begin{equation}\label{u5}
\psi(s)=\int\Big(\int\kappa(s)\n(s)ds\Big)ds+C,
\end{equation}
or in parametric representation
\begin{equation}\label{u6}
\psi(\theta)=\int\frac{1}{\kappa(\theta)}\Big(\int\n(\theta)d\theta\Big)d\theta+C,
\end{equation}
where $\theta=\int\kappa(s)ds$.

We can solve the equation (\ref{u2}) in the case of a slant helix. Firstly, we give:

\begin{lemma} \label{lm-1} Let $\psi:I\rightarrow \e^3$ be a curve that is parameterized by arclength with the intrinsic equations $\kappa=\kappa(s)$ and $\tau=\tau(s)$. The curve is a slant helix (its normal vectors make a constant angle, $\phi=\pm\mathrm{arccos}[n]$, with a fixed straight line in the space) if and only if
\begin{equation}\label{u7}
\tau(s)=\pm\dfrac{m\,\kappa(s)\int\kappa(s)ds}{\sqrt{1-m^2\Big(\int\kappa(s)ds\Big)^2}},
\end{equation}
where $m=\frac{n}{\sqrt{1-n^2}}$.
\end{lemma}

{\bf proof:} $(\Rightarrow)$ Let $\textbf{d}$ be the unitary fixed vector makes a constant angle, $\phi=\pm\mathrm{arccos}[n]$, with the normal vector $\n$. Therefore
\begin{equation}\label{u8}
\langle\n,\textbf{d}\rangle=n.
\end{equation}
Differentiating the equation (\ref{u8}) with respect to the variable $\theta=\int\kappa(s)ds$ and using the new Frenet equations (\ref{u3}), we get
\begin{equation}\label{u9}
\langle-\t(\theta)+f(\theta)\b(\theta),\textbf{d}\rangle=0.
\end{equation}
Therefore,
$$
\langle\t,\textbf{d}\rangle=f\langle\b,\textbf{d}\rangle.
$$
If we put $\langle\b,\textbf{d}\rangle=b$, we can write
$$
\textbf{d}=f\,b\,\t+n\,\n+b\,\b.
$$
From the unitary of the vector $\textbf{d}$ we get $b=\pm\sqrt{
\dfrac{1-n^2}{1+f^2}}$. Therefore, the vector $\textbf{d}$ can be written as
\begin{equation}\label{u10}
\textbf{d}=\pm\,f\,\sqrt{
\dfrac{1-n^2}{1+f^2}}\,\t+n\,\n\pm\sqrt{
\dfrac{1-n^2}{1+f^2}}\,\b.
\end{equation}
If we differentiate Equation (\ref{u5}) again, we obtain
\begin{equation}\label{u11}
\langle f'\,\b+(1+f^2)\n,\textbf{d}\rangle=0.
\end{equation}
Equations (\ref{u10}) and (\ref{u11}) lead to the following differential equation
$$
\frac{f'}{(1+f^2)^{3/2}}=\pm\,m,
$$
where $m=\frac{n}{\sqrt{1-n^2}}$. Integration the above equation, we get
\begin{equation}\label{u12}
\frac{f}{\sqrt{1+f^2}}=\pm\,m\,(\theta+c_1).
\end{equation}
where $c_1$ is an integration constant. The integration constant can disappear with a parameter change $\theta\rightarrow\theta-c_1$. Solving the equation (\ref{u12}) with $f$ as unknown we have
\begin{equation}\label{u13}
f(\theta)=\pm\frac{m\,\theta}{\sqrt{1-m^2\theta^2}}.
\end{equation}
Finally, $\tau(s)=\kappa(s)f(s)$, we express the desired result.

$(\Leftarrow)$ Suppose that $\tau(s)=\pm\dfrac{m\,\kappa(s)\int\kappa(s)ds}{\sqrt{1-m^2\Big(\int\kappa(s)ds\Big)^2}}$. The function $f$ can be written as $f(\theta)=\pm\frac{m\,\theta}{\sqrt{1-m^2\theta^2}}$ and let us consider the vector
$$
\textbf{d}=n\Big(\theta\,\t+\n\pm\frac{1}{m}\sqrt{1-m^2\theta^2}\,\b\Big).
$$
We will prove that the vector $\textbf{d}$ is a constant vector. Indeed, applying Frenet formula (\ref{u3})
$$
\textbf{d}'=n\Big(\t+\theta\,\n-\t+f\b\mp\frac{m\,\theta}{\sqrt{1-m^2\theta^2}}\b
\mp\frac{f}{m}\sqrt{1-m^2\theta^2}\n\Big)=0
$$
Therefore, the vector $\textbf{d}$ is constant and $\langle\n,\textbf{d}\rangle=n$. This concludes the proof of Lemma (\ref{lm-1}).

\begin{theorem}\label{th-main2} The position vector $\psi=(\psi_1,\psi_2,\psi_3)$ of a slant helix is computed in the natural representation form:
\begin{equation}\label{u14}
\left\{
\begin{array}{ll}
\psi_1(s)=\frac{n}{m}\,\int\Bigg[\int\kappa(s)\cos\Big[\frac{1}{n}
          \mathrm{arcsin}\Big(m\int\kappa(s)ds\Big)\Big]ds\Bigg]ds,\\
\psi_2(s)=\frac{n}{m}\,\int\Bigg[\int\kappa(s)\sin\Big[\frac{1}{n}
          \mathrm{arcsin}\Big(m\int\kappa(s)ds\Big)\Big]ds\Bigg]ds,\\
\psi_3(s)=n\,\int\Big[\int\kappa(s)ds\Big]ds,
\end{array}
\right.
\end{equation}
or in the parametric form
\begin{equation}\label{u15}
\left\{
\begin{array}{ll}
\psi_1(\theta)=\frac{n}{m}\,\int\frac{1}{\kappa(\theta)}\Bigg[\int\cos\Big[\frac{1}{n}
          \mathrm{arcsin}\Big(m\theta\Big)\Big]d\theta\Bigg]d\theta,\\
\psi_2(\theta)=\frac{n}{m}\,\int\frac{1}{\kappa(\theta)}\Bigg[\int\sin\Big[\frac{1}{n}
          \mathrm{arcsin}\Big(m\theta\Big)\Big]d\theta\Bigg]d\theta,\\
\psi_3(\theta)=n\,\int\frac{\theta}{\kappa(\theta)}d\theta,
\end{array}
\right.
\end{equation}
or in the useful parametric form:
\begin{equation}\label{u16}
\left\{
\begin{array}{ll}
\psi_1(t)=\frac{n^3}{m^3}\,\int\frac{\cos[nt]}{\kappa(t)}\Big[\int\cos[t]\cos[nt]dt\Big]dt,\\
\psi_2(t)=\frac{n^3}{m^3}\,\int\frac{\cos[nt]}{\kappa(t)}\Big[\int\sin[t]\cos[nt]dt\Big]dt,\\
\psi_3(t)=\frac{n^2}{m^2}\,\int\frac{\sin[nt]\cos[nt]}{\kappa(t)}dt,
\end{array}
\right.
\end{equation}
where $\theta=\int\kappa(s)ds$, $t=\frac{1}{n}\mathrm{arcsin}(m\theta)$, $m=\frac{n}{\sqrt{1-n^2}}$, $n=\cos[\phi]$ and $\phi$ is the angle between the fixed straight line (axis of a slant helix) and the principal normal vector of the curve.
\end{theorem}

{\bf Proof:} If $\psi$ is a slant helix whose principal normal vector $\n$ makes an angle $\phi=\pm\mathrm{\arccos}[n]$ with a straight line $U$, then we can write $f(\theta)=\pm
\dfrac{m\,\theta}{\sqrt{1-m^2\,\theta^2}}$, where $f=\frac{\tau}{\kappa}$, $\theta=\int\kappa(s)ds$ and $m=\frac{n}{\sqrt{1-n^2}}$. Therefore the equation (\ref{u2}) becomes
\begin{equation}\label{u17}
(1-m^2\theta^2)\n'''(\theta)-3m^2\theta\n''(\theta)+\n'(\theta)=0.
\end{equation}
If we write the principal normal vector as the following:
\begin{equation}\label{u18}
\n=N_1(\theta)\mathbf{e}_1+N_2(\theta)\mathbf{e}_2+N_3(\theta)\mathbf{e}_3.
\end{equation}
Now, the curve $\psi$ is a slant helix, i.e. the principal normal vector $\n$ makes a constant angle, $\phi$, with the constant vector called the axis of the slant helix. So, with out loss of generality, we take the axis of a slant helix is parallel to $\bold{e}_3$. Then
\begin{equation}\label{u19}
N_3=\langle\n,\mathbf{e}_3\rangle=n.
\end{equation}
On other hand the principal normal vector $\n$ is a unit vector, so the following condition is satisfied
\begin{equation}\label{u20}
N_1^2(\theta)+N_2^2(\theta)=1-n^2=\frac{n^2}{m^2}.
\end{equation}
The general solution of equation (\ref{u20}) can be written in the following form:
\begin{equation}\label{u21}
\begin{array}{ll}
N_1=\frac{n}{m}\cos[t(\theta)],\,\,\,\,\,N_2=\frac{n}{m}\sin[t(\theta)],
\end{array}
\end{equation}
where $t$ is an arbitrary function of $\theta$. Every component of the vector $\n$ is satisfied the equation (\ref{u17}). So, substituting the components $N_1(\theta)$ and $N_2(\theta)$ in the equation (\ref{u17}), we have the following differential equations of the function $t(\theta)$
\begin{equation}\label{u22}
3t'\Big(m^2\theta\,t'-(1-m^2\theta^2)t''\Big)\sin[t]-\Big(t'-3m^2\theta\,t''-(1-m^2\theta^2)(t'^3-t''')\Big)\cos[t]=0,
\end{equation}
\begin{equation}\label{u23}
3t'\Big(m^2\theta\,t'-(1-m^2\theta^2)t''\Big)\cos[t]+\Big(t'-3m^2\theta\,t''-(1-m^2\theta^2)(t'^3-t''')\Big)\sin[t]=0.
\end{equation}
It is easy to prove that the above two equations lead to the following two equations:
\begin{equation}\label{u24}
m^2\theta\,t'-(1-m^2\theta^2)t''=0,
\end{equation}
\begin{equation}\label{u25}
t'-3m^2\theta\,t''-(1-m^2\theta^2)(t'^3-t''')=0.
\end{equation}
The general solution of the equation (\ref{u24}) is
\begin{equation}\label{u26}
t(\theta)=c_2+c_1\,\mathrm{arcsin}(m\,\theta),
\end{equation}
or
\begin{equation}\label{u261}
t(\theta)=c_2+c_1\,\mathrm{arccos}(m\,\theta),
\end{equation}
where $c_1$ and $c_2$ are constants of integration. The constant $c_2$ can be disappear if we change the parameter $t\rightarrow t+c_2$. Substituting the solution (\ref{u26}) or (\ref{u261}) in the equation (\ref{u25}), we obtain the following condition:
$$
m\,c_1\,\Big(1+m^2(1-c_1)\Big)=0
$$
which leads to $c_1=\frac{\sqrt{1+m^2}}{m}=\frac{1}{n}$, where $m\neq 0$ and $c_1\neq 0$.

Now, the principal normal vector take the following form:
\begin{equation}\label{u27}
\n(\theta)=\Big(\frac{n}{m}\cos\Big[\frac{1}{n}\mathrm{arcsin}(m\,\theta)\Big],
\frac{n}{m}\sin\Big[\frac{1}{n}\mathrm{arcsin}(m\,\theta)\Big],n\Big).
\end{equation}
or
\begin{equation}\label{u271}
\n(\theta)=\Big(\frac{n}{m}\cos\Big[\frac{1}{n}\mathrm{arccos}(m\,\theta)\Big],
\frac{n}{m}\sin\Big[\frac{1}{n}\mathrm{arccos}(m\,\theta)\Big],n\Big).
\end{equation}
If we substitute the equation (\ref{u27}) in the two equations (\ref{u5}) and (\ref{u6}), we have the two equations (\ref{u14}) and (\ref{u15}). It is easy to arrive the equation (\ref{u16}), if we take the new parameter $t=\frac{1}{n}\mathrm{arcsin}(m\theta)$, which it completes the proof.

On other hand if we used the equation (\ref{u271}), we have the following theorem:

\begin{theorem}\label{th-main3} The position vector $\psi=(\psi_1,\psi_2,\psi_3)$ of a slant helix is computed in the natural representation form:
\begin{equation}\label{u141}
\left\{
\begin{array}{ll}
\psi_1(s)=\frac{n}{m}\,\int\Bigg[\int\kappa(s)\cos\Big[\frac{1}{n}
          \mathrm{arccos}\Big(m\int\kappa(s)ds\Big)\Big]ds\Bigg]ds,\\
\psi_2(s)=\frac{n}{m}\,\int\Bigg[\int\kappa(s)\sin\Big[\frac{1}{n}
          \mathrm{arccos}\Big(m\int\kappa(s)ds\Big)\Big]ds\Bigg]ds,\\
\psi_3(s)=n\,\int\Big[\int\kappa(s)ds\Big]ds,
\end{array}
\right.
\end{equation}
or in the parametric form
\begin{equation}\label{u151}
\left\{
\begin{array}{ll}
\psi_1(\theta)=\frac{n}{m}\,\int\frac{1}{\kappa(\theta)}\Bigg[\int\cos\Big[\frac{1}{n}
          \mathrm{arccos}\Big(m\theta\Big)\Big]d\theta\Bigg]d\theta,\\
\psi_2(\theta)=\frac{n}{m}\,\int\frac{1}{\kappa(\theta)}\Bigg[\int\sin\Big[\frac{1}{n}
          \mathrm{arccos}\Big(m\theta\Big)\Big]d\theta\Bigg]d\theta,\\
\psi_3(\theta)=n\,\int\frac{\theta}{\kappa(\theta)}d\theta,
\end{array}
\right.
\end{equation}
or in the useful parametric form:
\begin{equation}\label{u161}
\left\{
\begin{array}{ll}
\psi_1(t)=\frac{n^3}{m^3}\,\int\frac{\sin[nt]}{\kappa(t)}\Big[\int\cos[t]\sin[nt]dt\Big]dt,\\
\psi_2(t)=\frac{n^3}{m^3}\,\int\frac{\sin[nt]}{\kappa(t)}\Big[\int\sin[t]\sin[nt]dt\Big]dt,\\
\psi_3(t)=-\frac{n^2}{m^2}\,\int\frac{\sin[nt]\cos[nt]}{\kappa(t)}dt,
\end{array}
\right.
\end{equation}
where $\theta=\int\kappa(s)ds$, $t=\frac{1}{n}\mathrm{arccos}(m\theta)$, $m=\frac{n}{\sqrt{1-n^2}}$, $n=\cos[\phi]$ and $\phi$ is the angle between the fixed straight line (axis of a slant helix) and the principal normal vector of the curve.
\end{theorem}

\section{Examples}

In this section, we take several choices for the curvature $\kappa$ and torsion $\tau$, and next, we apply Theorem \ref{th-main2}.

\begin{example}
The case of a slant helix with
\begin{equation}\label{u28}
\kappa=1,\,\,\,\,\,\tau=\frac{m\,s}{\sqrt{1-m^2\,s^2}},
\end{equation}
which is the intrinsic equations of a Salkowski curve \cite{mont1}. Substituting $\kappa(t)=1$ in the equation (\ref{u16}) we have the explicit parametric representation of such curve as follows:
\begin{equation}\label{u281}
\left\{
\begin{array}{ll}
\psi_1(t)=\frac{n}{4m}\Big[\frac{n-1}{2n+1}\cos[(2n+1)t]+
\frac{n+1}{2n-1}\cos[(2n-1)t]-2\cos[t]\Big],\\
\psi_2(t)=\frac{n}{4m}\Big[\frac{n-1}{2n+1}\sin[(2n+1)t]-
\frac{n+1}{2n-1}\sin[(2n-1)t]-2\sin[t]\Big],\\
\psi_3(t)=-\frac{n}{4m^2}\cos[2nt].
\end{array}
\right.
\end{equation}
The curvature of the above curve is equal $1$ and the torsion is equal
$$
\tau=\tan[nt]=\frac{m\,\theta}{\sqrt{1-m^2\,\theta^2}}=\frac{m\,s}{\sqrt{1-m^2\,s^2}}.
$$
where $\theta=s$ and $t=\frac{1}{n}\mathrm{arcsin}(m\theta)$. One can see a special examples of such curves in the figure 1.
\begin{figure}[ht]
\begin{center}
\includegraphics[width=4.2cm]{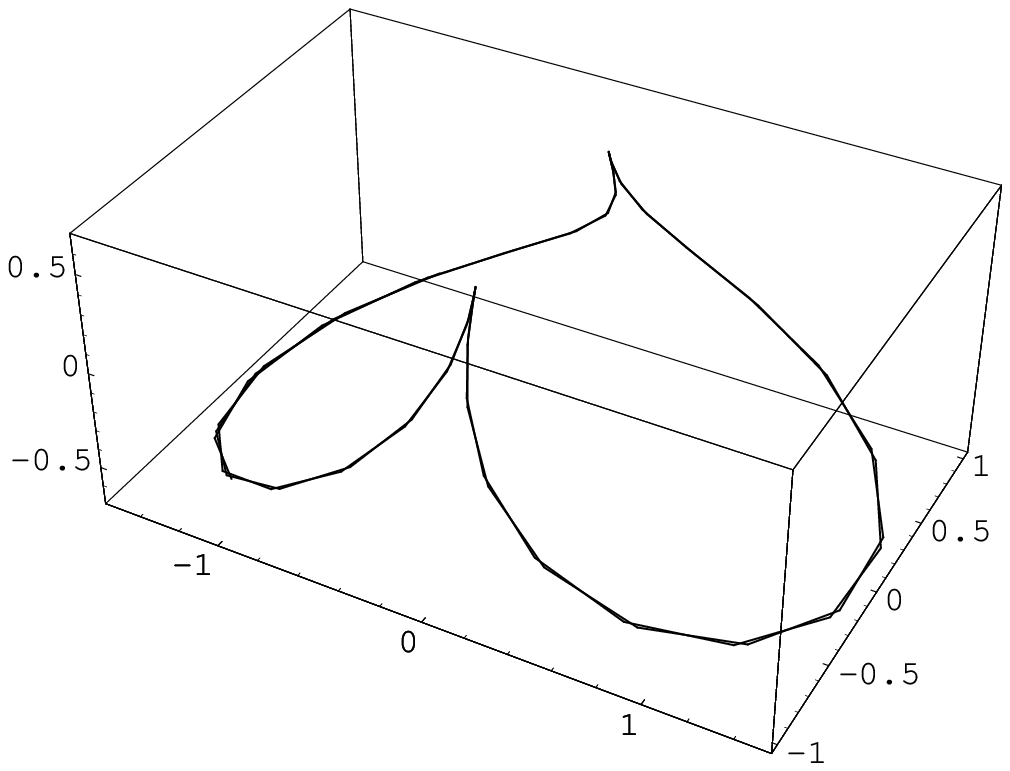}
\hspace*{0.5cm}
\includegraphics[width=3.2cm]{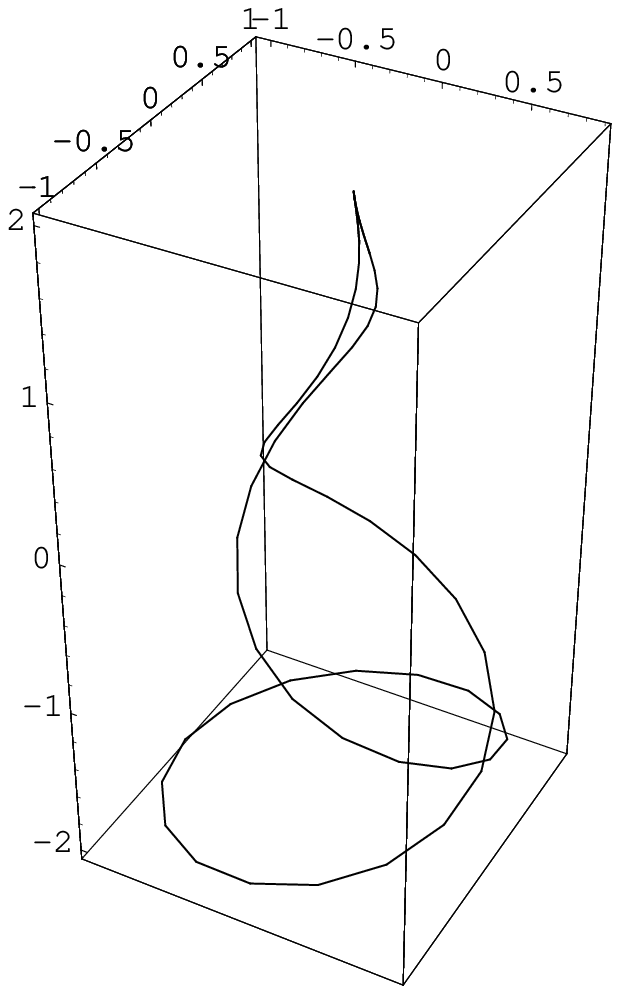}
\hspace*{0.5cm}
\includegraphics[width=5.2cm]{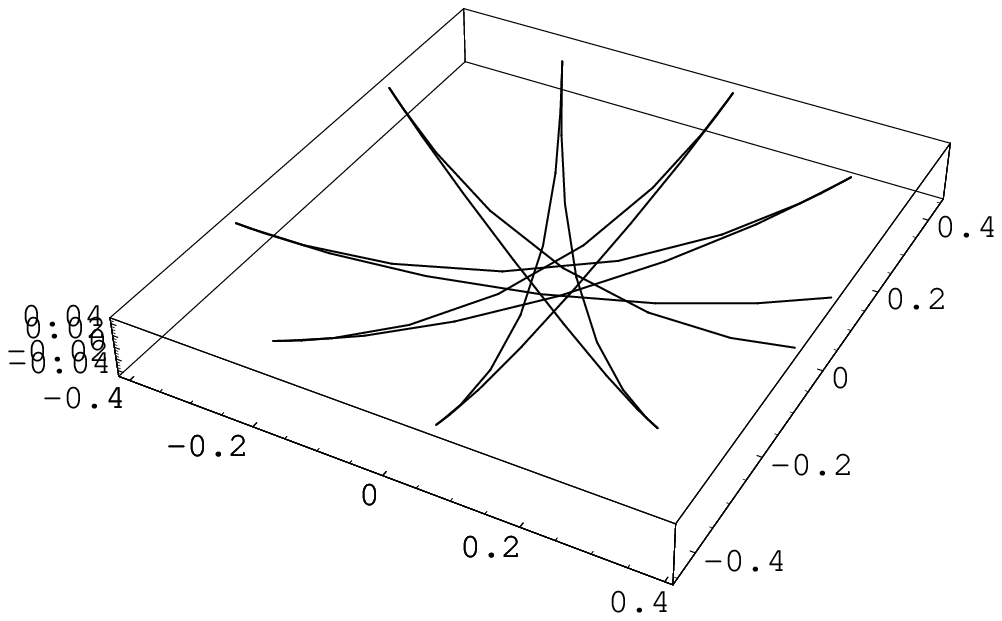}
\end{center}\caption{Some Slant helices with $\kappa=1$ and $n=\frac{1}{3}, \frac{1}{8}, \frac{10}{11}$.}
\label{fig-x1}
\end{figure}
\end{example}

\begin{example}
The case of a slant helix with
\begin{equation}\label{u29}
\kappa=\frac{m\,s}{\sqrt{1-m^2\,s^2}},\,\,\,\,\,\tau=1.
\end{equation}
which is the intrinsic equations of an anti-Salkowski curve \cite{mont1}. Substituting
$$
\kappa=\frac{m\,s}{\sqrt{1-m^2\,s^2}}=\frac{\sqrt{1-m^2\theta^2}}{m\theta}=\cot[nt],
$$
in the equation (\ref{u16}) we have the explicit parametric representation of such curve as follows:
\begin{equation}\label{u30}
\left\{
\begin{array}{ll}
\psi_1(t)=\frac{n}{4m}\Big[\frac{n-1}{2n+1}\sin[(2n+1)t]+
\frac{n+1}{2n-1}\sin[(2n-1)t]-2n\sin[t]\Big],\\
\psi_2(t)=\frac{n}{4m}\Big[\frac{1-n}{1+2n}\cos[(1+2n)t]-
\frac{1+n}{1-2n}\cos[(1-2n)t]+2n\cos[t]\Big],\\
\psi_3(t)=\frac{n}{4m^2}(2nt-\sin[2nt]),
\end{array}
\right.
\end{equation}
where $\theta=\frac{\sqrt{1-m^2s^2}}{m}$ and $t=\frac{1}{n}\mathrm{arcsin}(m\theta)$. One can see a special examples of such curves in the figure 2.
\begin{figure}[ht]
\begin{center}
\includegraphics[width=3.2cm]{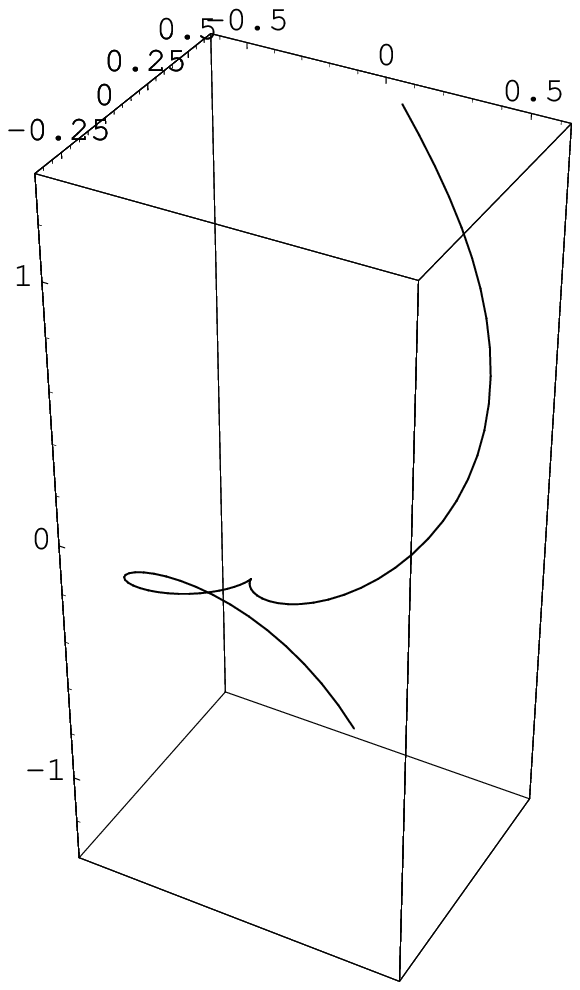}
\hspace*{0.5cm}
\includegraphics[width=4.2cm]{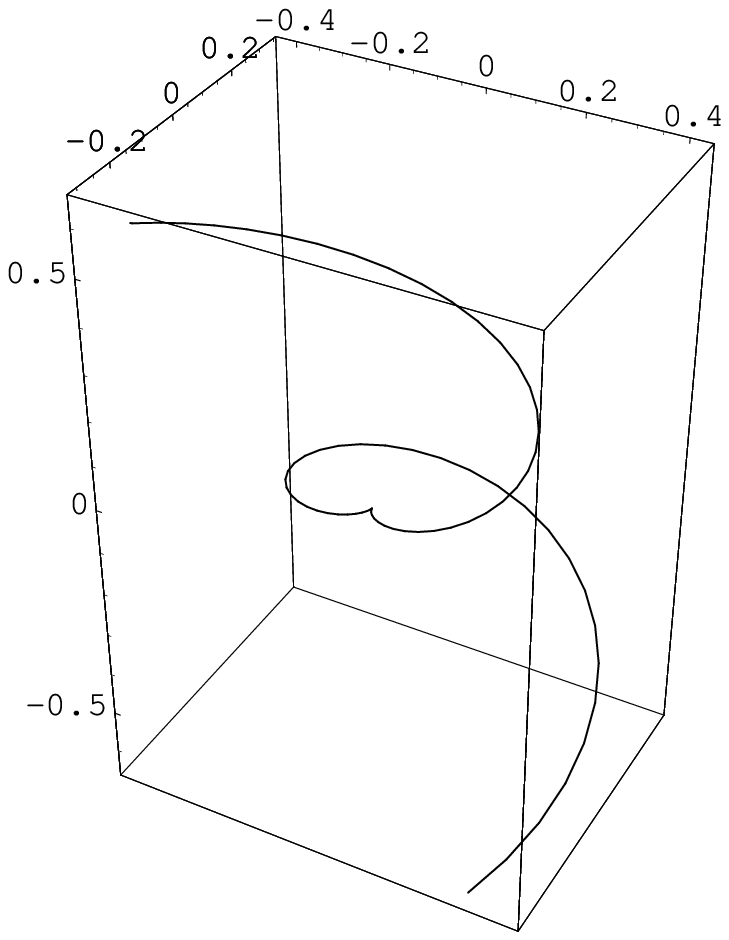}
\hspace*{0.5cm}
\includegraphics[width=3.2cm]{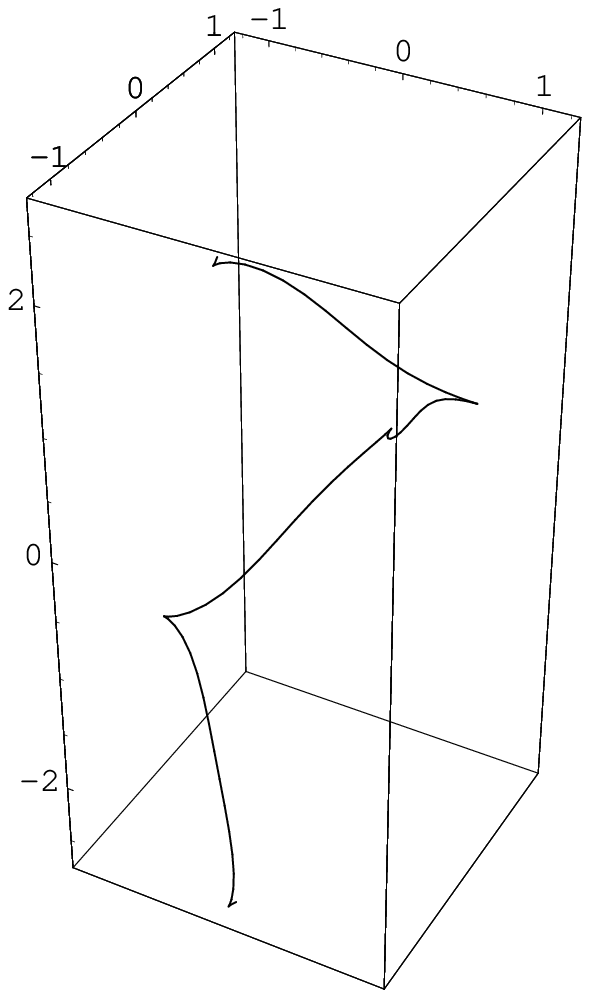}
\end{center}\caption{Some Slant helices with $\tau=1$ and $n=\frac{1}{5}, \frac{1}{13}, \frac{2}{3}$.}
\label{fig-x2}
\end{figure}
\end{example}

\begin{remark}
 A family of curves with constant curvature but non-constant torsion is called Salkowski curves and a family of curves with constant torsion but non-constant curvature is called anti-Salkowski curves $\cite{salkow}$. Monterde \cite{mont1} studied some of characterizations of these curves and he prove that the principal normal vector makes a constant angle with fixed straight line. So that: Salkowski and anti-Salkowski curves are the important examples of slant helices.
\end{remark}

\begin{example}
The case of a slant helix with
\begin{equation}\label{u31}
\kappa=\frac{\mu}{m}\cos[\mu\,s],\,\,\,\,\,\tau=\frac{\mu}{m}\sin[\mu\,s].
\end{equation}
Substituting $\kappa=\frac{\mu}{m}\cos[m\,s]$ in the equation (\ref{u14}) we have the natural representation of such curve as follows:
\begin{equation}\label{u32}
\left\{
\begin{array}{ll}
\psi_1(s)=-\frac{m^2}{n\,\mu}\Big[(1+n^2)\cos[\mu\,s]\cos[\frac{\mu\,s}{n}]+
2n\sin[\mu\,s]\sin[\frac{\mu\,s}{n}]\Big],\\
\psi_2(s)=-\frac{m^2}{n\,\mu}\Big[(1+n^2)\cos[\mu\,s]\sin[\frac{\mu\,s}{n}]-
2n\sin[\mu\,s]\cos[\frac{\mu\,s}{n}]\Big],\\
\psi_3(s)=-\frac{n}{m\,\mu}\,\cos[\mu\,s].
\end{array}
\right.
\end{equation}
The above curve is a geodesic of the tangent developable of a general helix \cite{izumi}. One can see a special examples of such curves in the figure 3.
\begin{figure}[ht]
\begin{center}
\includegraphics[width=4.7cm]{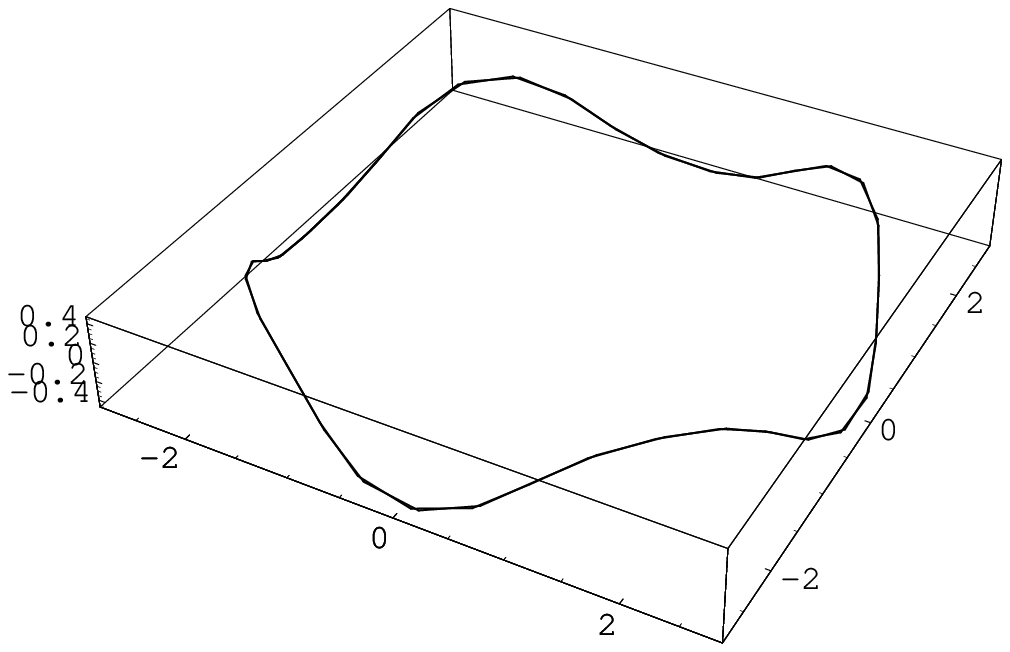}
\hspace*{0.5cm}
\includegraphics[width=3.7cm]{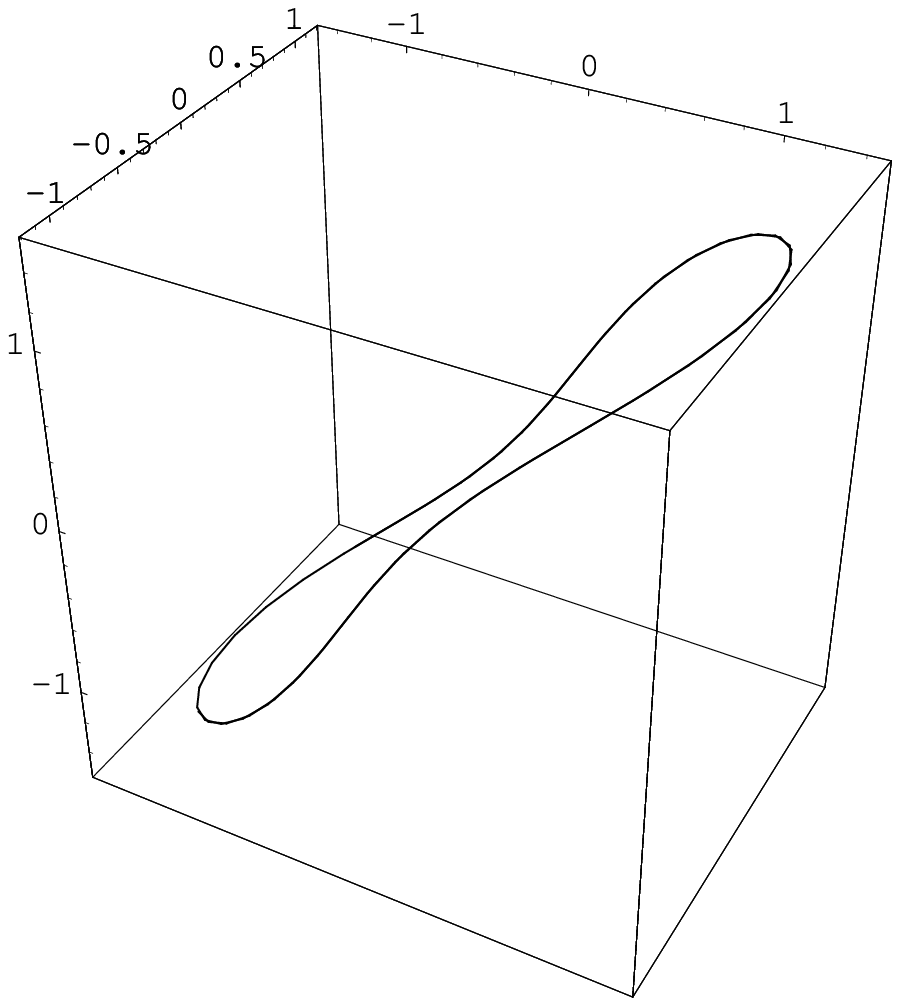}
\hspace*{0.5cm}
\includegraphics[width=2.7cm]{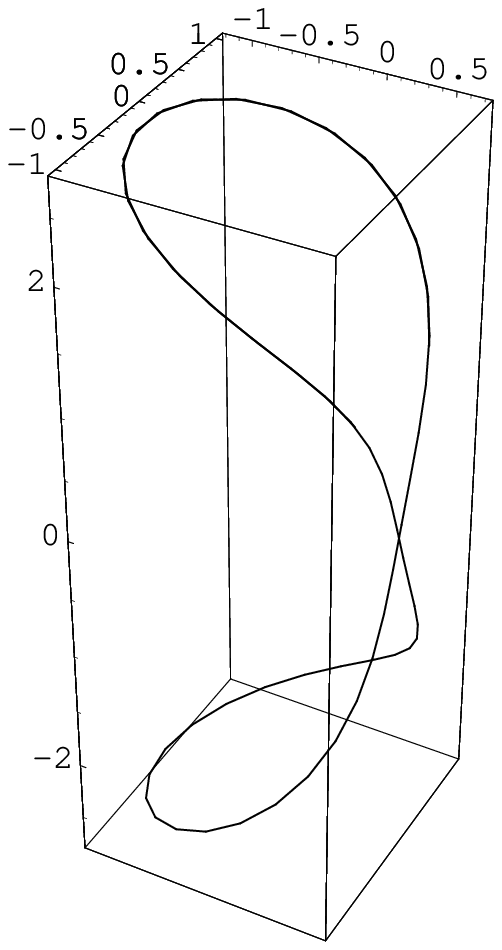}
\end{center}\caption{Some Slant helices with $\kappa=\frac{\mu}{m}\cos[\mu\,s]$, $\mu=m$ and $n=\frac{4}{5}, \frac{1}{2}, \frac{1}{3}$.}
\label{fig-x3}
\end{figure}
\end{example}

\begin{remark}
A unit speed curve of constant precession is defined by the property that its (Frenet) Darboux vector
$$
W=\tau\,\t+\kappa\,\b
$$
revolves about a fixed line in space with angle and constant speed. A curve of constant precession is characterized by having
$$
\kappa=\frac{\mu}{m}\sin[\mu\,s],\,\,\,\,\,\,\,\,\,\,\tau=\frac{\mu}{m}\cos[\mu\,s]
$$
or
$$
\kappa=\frac{\mu}{m}\cos[\mu\,s],\,\,\,\,\,\,\,\,\,\,\tau=\frac{\mu}{m}\sin[\mu\,s]
$$
where $\mu$ and $m$ are constants. This curve lie on a circular one-sheeted hyperboloid
$$
x^2+y^2-m^2\,z^2=4m^2.
$$
The curve closed if and only if $n=\frac{m}{\sqrt{1+m^2}}$ is rational \cite{sc}. Kula and Yayli \cite{kula} proved that the geodesic curvature of the spherical image of the principal normal indicatrix of a curve of constant precession is a constant function equal $-m$. So that: the curves of constant precessions are the important examples of slant helices.
\end{remark}


\end{document}